\documentclass[12pt,reqno]{amsart}

\usepackage{amssymb}
\usepackage{amsmath}

\newcommand{\alp}{\alpha}
\newcommand{\del}{\delta}
\newcommand{\kap}{\kappa}
\renewcommand{\phi}{\varphi}
\newcommand{\sig}{\sigma}

\newcommand{\Gam}{\Gamma}

\newcommand{\F}{{\mathbb F}}

\newcommand{\Fr}[1][r]{\F_2^{#1}}

\newcommand{\tD}{{\widetilde D}}

\renewcommand{\(}{\left(}
\renewcommand{\)}{\right)}

\newcommand{\rfr}{\right\}}

\newcommand{\lcl}{\left\lceil}
\newcommand{\rcl}{\right\rceil}

\renewcommand{\>}{\right\rangle}

\newcommand{\seq}{\subseteq}
\newcommand{\stm}{\setminus}
\newcommand{\est}{\varnothing}

\newcommand{\longc}{,\dotsc,}

\newcommand{\multplus}[1][2]{\stackrel{#1}{+}}

\DeclareMathOperator{\PG}{PG}

\newtheorem{theorem}{Theorem}
\newtheorem{fact}[theorem]{Fact}
\newtheorem{claim}[theorem]{Claim}
\newtheorem{lemma}[theorem]{Lemma}
\newtheorem{corollary}[theorem]{Corollary}
\newtheorem{proposition}[theorem]{Proposition}

\newtheorem{primetheorem}{Theorem}
\newtheorem{alphatheorem}{Theorem}

\theoremstyle{remark}
\newtheorem*{remark}{Remark}

\theoremstyle{definition}
\newtheorem{example}[theorem]{Example}

\newcommand{\reff}[1]{\ref{f:#1}}
\newcommand{\refc}[1]{\ref{c:#1}}
\newcommand{\reft}[1]{\ref{t:#1}}
\newcommand{\refl}[1]{\ref{l:#1}}
\newcommand{\refm}[1]{\ref{m:#1}}
\newcommand{\refp}[1]{\ref{p:#1}}
\newcommand{\refx}[1]{\ref{x:#1}}
\newcommand{\refs}[1]{\ref{s:#1}}
\newcommand{\refb}[1]{\cite{b:#1}}
\newcommand{\refe}[1]{\eqref{e:#1}}

\author{David J. Grynkiewicz}
\email{diambri@hotmail.com}
\address{Institut f\"ur Mathematik und Wissenschaftliches Rechnen \\
         Karl-Franzens-Universit\"at Graz \\
         Heinrichstra\ss e 36 \\
         8010 Graz, Austria}
\thanks{The first author was supported by the Austrian Science Fund FWF
       (Project Number M1014-N13)}

\author{Vsevolod F. Lev}
\email{seva@math.haifa.ac.il}
\address{Department of Mathematics, The University of Haifa at Oranim,
         Tivon 36006, Israel}

\title[$1$-Saturating Sets, Caps and Round Sets]%
  {$1$-Saturating Sets, Caps and Round Sets \\ in Binary Spaces}


\begin{document}
\baselineskip=16pt

\begin{abstract}
We show that, for a positive integer $r$, every minimal $1$-saturating set in
$\PG(r-1,2)$ of size at least $\frac{11}{36}\,2^r+3$ either is a complete cap
or can be obtained from a complete cap $S$ by fixing some $s\in S$ and
replacing every point $s'\in S\stm\{s\}$ by the third point on the line
through $s$ and $s'$. Stated algebraically: if $G$ is an elementary abelian
$2$-group and a set $A\seq G\stm\{0\}$ with $|A|>\frac{11}{36}\,|G|+3$
satisfies $A\cup 2A=G$ and is minimal subject to this condition, then either
$A$ is a maximal sum-free set, or there are a maximal sum-free set $S\seq G$
and an element $s\in S$ such that $A=\{s\}\cup\big(s+(S\stm\{s\})\big)$.
Since, conversely, every set obtained in this way is a minimal $1$-saturating
set, and the structure of large sum-free sets in an elementary abelian
$2$-group is known, this provides a complete description of large minimal
$1$-saturating sets.

Our approach is based on characterizing those large sets $A$ in elementary
abelian $2$-groups such that, for every proper subset $B$ of $A$, the sumset
$2B$ is a proper subset of $2A$.
\end{abstract}

\maketitle

\section{Saturating Sets and Caps: The Main Result.}\label{s:intro}

Let $r\ge 1$ be an integer, $q$ a prime power, and $A\seq\PG(r-1,q)$ a set in
the $(r-1)$-dimensional projective space over the $q$-element field. Given an
integer $\rho\ge 1$, one says that $A$ is \emph{$\rho$-saturating} if every
point of $\PG(r-1,q)$ is contained in a subspace generated by $\rho+1$ points
from $A$. Furthermore, $A$ is said to be a \emph{cap} if no three points of
$A$ are collinear; a cap is \emph{complete} if it is not properly contained
in another cap. Since the property of being $\rho$-saturating is inherited by
supersets  and that of being a cap is inherited by subsets, of particular
interest are minimal $\rho$-saturating sets and complete caps.

In this paper, we are concerned with the case $\rho=1$ and the space
$\PG(r-1,2)$ whose points are, essentially, non-zero elements of the
elementary abelian $2$-group of rank $r$, and whose lines are triples of
points adding up to $0$. A large random set in $\PG(r-1,2)$ is $1$-saturating
with very high probability, but the probability that it is a \emph{minimal}
$1$-saturating set is extremely low; thus, one can expect that large minimal
$1$-saturating sets are well-structured and can be explicitly described. A
similar heuristic applies to large complete caps, and indeed, a classical
result of Davydov and Tombak \refb{dt} establishes the structure of complete
caps of size larger than $2^{r-2}+1$. Classifying large $1$-saturating sets
seems to be considerably more subtle, which is
quite natural bearing in mind that complete caps in $\PG(r-1,2)$ can be
characterized as those $1$-saturating sets possessing the
extra property of having no internal lines (as will be explained shortly).

With the exception of the next section where our result is discussed from the
projective geometric viewpoint, we mostly use the language of abelian groups.
Accordingly, denoting by $\Fr$ the elementary abelian $2$-group of rank $r\ge
1$ and writing
  $$ 2A := \{a_1+a_2\colon a_1,a_2\in A\} $$
for a subset $A\seq\Fr$, we interpret $1$-saturating sets in $\PG(r-1,2)$ as
those subsets $A\seq\Fr\stm\{0\}$
satisfying $A\cup 2A=\Fr$. Similarly, caps in $\PG(r-1,2)$ are understood as
sets $A\seq\Fr\stm\{0\}$ with $A\cap 2A=\est$; such sets are customarily
referred to as \emph{sum-free}. Complete caps are thus identified with
maximal (by inclusion) sum-free sets.

It is well known and easy to see that a sum-free set $A\seq\Fr$ is maximal if
and only if the sets $A$ and $2A$ partition $\Fr$; that is, in addition to
being disjoint, they satisfy $A\cup 2A=\Fr$. Consequently, any maximal
sum-free set is a minimal $1$-saturating set without internal lines. Beyond
this simple observation, the only general result which seems to be known
about minimal $1$-saturating sets in $\Fr$ is established in \refb{dmp1}; it
asserts that the largest possible size of such a set is $2^{r-1}$, examples
being furnished by the following two constructions:
\begin{itemize}
\item[(i)]  if $H<\Fr$ is an index-$2$ subgroup and $g\in\Fr\stm H$, then
    $g+H$ is a minimal $1$-saturating set;
\item[(ii)] with $H$ and $g$ as in (i), the union $\{g\}\cup(H\stm\{0\})$
    is a minimal $1$-saturating set.
\end{itemize}

An extension of construction (i) has just been mentioned: any maximal
sum-free set is a minimal $1$-saturating set. Construction (ii) can be
extended by observing that if $S$ is a maximal sum-free set and $s\in S$,
then $A:=\{s\}\cup\big((S+s)\stm\{0\}\big)$ is a minimal $1$-saturating set:
for in this case,
  $$ A\cup 2A = 2(A\cup\{0\}) = 2\big(s + (S\cup\{0\})\big)
                                      = 2(S\cup\{0\}) = S\cup 2S = \Fr, $$
and this computation also shows that, for any proper subset $B\subset A$, we
have $B\cup2B\neq\Fr$.

Indeed, a common description can be given to these two extensions: namely, if
$S\seq\Fr$ is a maximal sum-free set and $s\in S\cup\{0\}$, then
$A:=\big(s+(S\cup\{0\})\big)\stm\{0\}$ is a minimal $1$-saturating set. In
this paper, we classify completely minimal $1$-saturating sets in $\Fr$ of
size at least $\frac{11}{36}\,2^r+3$, showing that they all are of this form.

%\newcounter{MainSaturatingTheorem}
%\setcounter{MainSaturatingTheorem}{\value{theorem}}
\begin{theorem}\label{t:main}
Let $r\ge 1$ be an integer. A set $A\seq\Fr\stm\{0\}$ with
$|A|>\frac{11}{36}\,2^r+3$ is a minimal $1$-saturating set if and only if
there are a maximal sum-free set $S\seq\Fr$ and an element $s\in S\cup\{0\}$
such that $A=\big(s+(S\cup\{0\})\big)\stm\{0\}$.
\end{theorem}

We notice that Theorem \reft{main} provides a comprehensive characterization
of large minimal $1$-saturating sets, as the structure of large maximal
sum-free sets is known due to the result of Davydov and Tombak mentioned at
the beginning of this section. We record the following easy corollary of
their result.

\begin{fact}[\cite{b:dt}]\label{f:dt}
Let $r\ge 1$ be an integer. Every maximal sum-free set in $\Fr$ of size
larger than $9\cdot 2^{r-5}$ either is the non-zero coset of an index-$2$
subgroup, or has the form $B+H$, where $H<\Fr$ is a subgroup of index $16$
and $B\subset\Fr$ is a five-element set with $\Fr=\<B\>\oplus H$ such that
the elements of $B$ add up to $0$.
\end{fact}
In the statement of Fact~\reff{dt} and below in the paper, for a set $B$ of
group elements, we use $\<B\>$ to denote the subgroup generated by $B$.
Furthermore, given yet another subset $C$ of the same group, we write
$B+C:=\{b+c\colon b\in B,\,c\in C\}$. The set $B+C$ is commonly referred to
as the \emph{sumset} of $B$ and $C$. Notice that $B+B=2B$.

We conjecture that the density assumption of Theorem \reft{main} can actually
be relaxed to $|A|\ge 2^{r-2}+3$, provided that $r\ge 6$. (The group $\F_2^5$
contains an $11$-element minimal $1$-saturating set, but no $11$-element
maximal sum-free sets; see \refb{dmp2}.) If true, this is best possible.

\begin{example}\label{x:triangle}
Given an integer $r\ge 4$, fix elements $e_1,e_2\in\Fr$ and an index-$4$
subgroup $H<\Fr$ with $\Fr=\<e_1,e_2\>\oplus H$, and let $A:=(\<e_1,e_2\>\cup
H)\stm\{0\}$. Straightforward verification shows that $A$ is a minimal
$1$-saturating set. Now, if $A=\big(s+(S\cup\{0\})\big)\stm\{0\}$ for a
subset $S\seq\Fr\stm\{0\}$ and an element $s\in S\cup\{0\}$, then
$S\cup\{0\}=s+(\<e_1,e_2\>\cup H)$. Since this set contains $0$, we have
$s\in\<e_1,e_2\>\cup H$. If $s\in H$, then $S$ contains all non-zero elements
of $H$, whence $2S=H$ in view of $|H|\ge 4$, and therefore $S$ is not
sum-free. If $s=e_1$, then $S=\{e_1,e_2,e_1+e_2\}\cup(e_1+H)$ is evidently
not sum-free, and similarly it is not sum-free if $s=e_2$ or $s=e_1+e_2$.
Thus $A$ cannot be represented as in Theorem~\reft{main}.

More generally, if $F$ and $H$ are subgroups with $\Fr=F\oplus H$ and
$|F|,|H|\ge 4$, then $(F\cup H)\stm\{0\}$ is a minimal $1$-saturating set
which cannot be represented as in Theorem~\reft{main}.
\end{example}

\section{The Projective Geometry Viewpoint}\label{s:PG}

We remark that Theorem \ref{t:main} can be reformulated in purely geometrical
terms, as in the abstract.
\renewcommand{\thealphatheorem}{\reft{main}a}
\begin{alphatheorem}\label{t:mainA}
For an integer $r\ge 1$, every minimal $1$-saturating set in $\PG(r-1,2)$ of
size at least $\frac{11}{36}\,2^r+3$ either is a complete cap, or can be
obtained from a complete cap $S$ by fixing some $s\in S$ and replacing every
point $s'\in S\stm\{s\}$ by the third point on the line through $s$ and $s'$.
\end{alphatheorem}

Another reformulation, kindly pointed out by Simeon Ball, involves blocking
sets. Recall that a set of points in a projective geometry is called a
\emph{blocking set} if it has a non-empty intersection with every line;
consequently, a set in $\PG(r-1,2)$ is a (minimal) blocking set if and only
if its complement is a (complete) cap. It is easy to derive that
Theorem~\reft{main} is equivalent to the following assertion.

\renewcommand{\thealphatheorem}{\reft{main}b}
\begin{alphatheorem}\label{t:mainB}
For an integer $r\ge 1$, every minimal $1$-saturating set $A$ in $\PG(r-1,2)$
of size at least $\frac{11}{36}\,2^r+3$ either is the complement of a minimal
blocking set, or can be obtained from a minimal blocking set $B$ by fixing a
point $s\notin B$ and letting $A$ consist of $s$ along with all points $b\in
B$ for which the line through $s$ and $b$ is tangent to $B$ (i.e., passes
through precisely one point of $B$).
\end{alphatheorem}
The equivalence between Theorems \reft{main} and \reft{mainB} is a simple
exercise, left to the reader.

Yet another consequence of Theorem~\reft{main} concerns the spectrum of
possible sizes of minimal $1$-saturating sets. As indicated in Section
\refs{intro}, the largest size of a minimal $1$-saturating set in
$\PG(r-1,2)$ is $2^{r-1}$. The \emph{second largest} size can be determined
as an immediate corollary of Theorem~\reft{main} and Fact~\reff{dt}.
\begin{corollary}
If $r\ge 9$ is an integer, then the second largest size of a minimal
$1$-saturating set in $\PG(r-1,2)$ is $5\cdot 2^{r-4}$, and the third largest
size is smaller than $\frac{11}{36}\,2^r+3$.
\end{corollary}

It is observed in \refb{dmp2} that, with a single exception for $r=5$, the
spectrum of sizes of all known large minimal $1$-saturating sets in
$\PG(r-1,2)$ is contained in the spectrum of sizes of sum-free sets in
$\PG(r-1,2)$. Theorem \reft{main} and its above-mentioned conjectured
strengthening provide, of course, an explanation to this phenomenon.

Finally, we note that Theorem \reft{main} allows one to find all classes of
projectively equivalent minimal $1$-saturating sets in $\PG(r-1,2)$. For, it
is not difficult to derive from Fact~\reff{dt} that if $r\ge 6$ is an
integer, $S_1$ and $S_2$ are (potentially identical) complete caps in
$\PG(r-1,2)$ with $|S_1|=|S_2|>9\cdot 2^{r-5}$, and, for $i\in\{1,2\}$, the
sets $S_i'$ are obtained from $S_i$ as described in Theorem \reft{mainA},
then $S_1$ and $S_2$ are projectively equivalent, as are $S_1'$ and $S_2'$,
while $S_1$ is not equivalent to $S_2'$---regardless of the specific choice
of the elements fixed in $S_1$ and $S_2$ to get $S_1'$ and $S_2'$ (for the
non-equivalence, one only needs to note that $S_2'$ is not a cap for
 $r\ge 6$.) This leads to the following corollary.
\begin{corollary}
For a positive integer $r\ge 9$, there are four projectively non-isomorphic
minimal $1$-saturating sets in $\PG(r-1,2)$ of size larger than
$\frac{11}{36}\,2^r+3$: two are complete caps of sizes $2^{r-1}$ and
$5\cdot2^{r-4}$, and two more are obtained from them as in Theorem
\reft{mainA}.
\end{corollary}

\section{Round Sets and the Unique Representation Graph.}\label{s:round}

In a paradoxical way, for a minimal $1$-saturating set, minimality seems to
be more important than saturation. This idea is captured in the notion of a
round set, introduced in the present section. We also bring into
consideration unique representation graphs, which are of fundamental
importance for our argument, and establish some basic properties of round
sets and unique representation graphs. Finally, we state a structure theorem
for round sets (Theorem~\reft{round} below) and show that it implies
Theorem~\reft{main}.

The remainder of the paper is structured as follows. Important auxiliary
results are gathered in Section \refs{aux}. In Section~\refs{light}, we prove
a ``light version'' of Theorem~\reft{main}, with the assumption on the size
of $A$ strengthened to $|A|>\frac13\,2^r+2$; besides supplying a proof of
Theorem~\reft{main} for small dimensions ($r\le 5$), it serves as a
simplified model of our method, exhibiting many of the core ideas.
Sections~\refs{twoisolated}--\refs{2edgescompletion} are devoted to the proof
of Theorem~\reft{round}: in Section~\refs{twoisolated}, the problem is
reduced to the case where the unique representation graph is known to have at
least two isolated edges, Sections \refs{2edgesbasics} and
\refs{2edgescompletion} present a treatment of this case.

Let $r\ge 1$ be an integer. We say that a set $A\seq\Fr$ is \emph{round} if,
for every proper subset $B\subset A$, we have $2B\ne 2A$; that is, for every
$a\in A$, there exists $a'\in A$ such that $a+a'$ has a unique (up to the
order of summands) representation as a sum of two elements of $A$.

It is immediate from the definition that $A\seq\Fr\stm\{0\}$ is a minimal
$1$-saturating set if and only if it satisfies $2(A\cup\{0\})=\Fr$ and is
minimal subject to this condition. The following simple lemma takes this
observation a little further.

\begin{lemma}\label{l:1S2R}
Let $r\ge 1$ be an integer. If $A\seq\Fr\stm\{0\}$ is a minimal
$1$-saturating set, then either $A$ or $A\cup\{0\}$ is round.
\end{lemma}

\begin{remark}
It is easy to derive from Theorem \reft{main} and the observation following
the proof below that if $A\seq\Fr\stm\{0\}$ is a large minimal $1$-saturating
set, then, indeed, $A\cup\{0\}$ is round.
\end{remark}

\begin{proof}[Proof of Lemma \refl{1S2R}]
Suppose that $A\seq\Fr$ is a minimal $1$-saturating set. If $A\cup\{0\}$ is
not round, then there exists $a_0\in A\cup\{0\}$ such that
$2\big((A\cup\{0\})\stm\{a_0\}\big)=2(A\cup\{0\})=\Fr$. Since $a_0\in A$
would contradict the minimality of $A$, we actually have $a_0=0$, whence
$2A=\Fr$. Now if also $A$ is not round, then there exists $a\in A$ with
$2(A\stm\{a\})=2A=\Fr$. This yields $2\big((A\stm\{a\})\cup\{0\}\big)=\Fr$,
which, again, contradicts the minimality of $A$.
\end{proof}

Lemma \refl{1S2R} allows us to concentrate on studying large round sets
instead of large $1$-saturating sets; indeed, we will hardly refer to
$1$-saturating sets from now on, except for the deduction of
Theorem~\reft{main} from Theorem~\reft{round} at the end of this section.

We observe that if $S\seq\Fr$ is sum-free, then $0\notin S$ and, for each
$g\in\Fr$, the set $g+(S\cup\{0\})$ is round. To verify this, we can assume
$g=0$ (as roundness is translation invariant) and notice that, fixing
arbitrarily $s_0\in S$ and letting $S_0:=S\stm\{s_0\}$, we have
 $s_0\notin 2S$ and $s_0\notin 2(S_0\cup\{0\})$, whereas
$s_0\in 2(S\cup\{0\})$. The heart of our paper is the following theorem,
showing that, in fact, any large round set has the structure just
described.
\newcounter{MainRoundTheorem}
\setcounter{MainRoundTheorem}{\value{theorem}}
\begin{theorem}\label{t:round}
Let $r\ge 1$ be an integer and suppose that $A\seq\Fr$ is round. If
$|A|>\frac{11}{36}\,2^r+3$, then there is a sum-free set $S\seq\Fr$ and
an element $g\in\Fr$ such that $A=g+(S\cup\{0\})$.
\end{theorem}

We now turn to the notion of a unique representation graph. Given an integer
$r\ge 1$ and a set $A\seq\Fr$, we define $D(A)$ to be the set of all those
elements of $\Fr$ with a unique, up to the order of summands, representation
as a sum of two elements of $A$.
By $\Gam(A)$ we denote the graph on the vertex set $A$ in which two vertices
$a_1,a_2\in A$ are adjacent whenever $a_1+a_2\in D(A)$; if $|A|>1$, then
$\Gam(A)$ is a simple, loopless graph (as all graphs below are tacitly
assumed to be). We call $\Gam(A)$ the \emph{unique representation graph} of
$A$. Notice that the number of edges of $\Gam(A)$ is $|D(A)|$ and that, for
any $g\in\Fr$, we have $D(A+g)=D(A)$, while $\Gam(g+A)$ is obtained from
$\Gam(A)$ by re-labeling the vertices.

Evidently, a set $A\seq\Fr$ with $|A|\ge 2$ is round if and only if $\Gam(A)$
has no isolated vertices. Another indication of the importance of unique
representation graphs is given by the following lemma.

\begin{lemma}\label{l:sf=star}
Let $r\ge 1$ be an integer, let $g\in\Fr$, and suppose that $A\seq\Fr$
satisfies $|A|\ge 2$. For $\Gam(A)$ to have a spanning star with the center
at $g$, it is necessary and sufficient that $A=g+(S\cup\{0\})$, where
$S\seq\Fr$ is sum-free.
\end{lemma}

\begin{proof}
If $g\notin A$, then $g$ is not a vertex of $\Gam(A)$ and $A\neq
g+(S\cup\{0\})$; thus, the assertion is immediate in this case. If $g\in A$,
set $S:=(A+g)\stm\{0\}$, so that $A=g+(S\cup\{0\})$. The graph $\Gam(A)$ has
a spanning star with the center at $g$ if and only if, for every $s\in S$, we
have $g+(g+s)\in D(A)$; that is, $g+(g+s)\ne (g+s_1)+(g+s_2)$ whenever
$s_1,s_2\in S$. This is equivalent to $S$ being sum-free.
\end{proof}

By Lemma \refl{sf=star}, to prove Theorem \reft{round}, it suffices to show
that if $A\seq\Fr$ is a large round set, then $\Gam(A)$ contains a spanning
star. The following basic result shows that, for the unique representation
graph of a large set, \emph{containing} a spanning star is equivalent to
\emph{being} a star.

\begin{proposition}\label{p:tfree}
Let $r\ge 1$ be an integer and suppose that $A\seq\Fr$. If
 $|A|\ge 2^{r-2}+3$, then $\Gam(A)$ is triangle-free. Moreover, if
$|A|>2^{r-2}+3$, then, indeed, $D(A)$ is sum-free.
\end{proposition}

\begin{remark}
Observe that if $a_1,a_2,a_3\in A$ induce a triangle in $\Gam(A)$, then
$D(A)$ is not sum-free in view of $(a_1+a_2)+(a_2+a_3)=a_1+a_3$; thus,
``$D(A)$ is sum-free'' is a stronger conclusion than ``$\Gam(A)$ is
triangle-free''. We also notice that the bound $2^{r-2}+3$ is sharp. To see
this, suppose that $e_1,\,e_2,\,H$, and $A$ are as in Example
\refx{triangle}, and set $A_0:=\<e_1,e_2\>\cup H$. Then $|A|=2^{r-2}+2$ and
the vertices $e_1,\,e_2$, and $e_1+e_2$ of $\Gam(A)$ induce a triangle,
whereas $|A_0|=2^{r-2}+3$ and $D(A_0)$ is not sum-free: for if $h_1$ and
$h_2$ are distinct non-zero elements of $H$, then $e_1+h_1,\,e_2+h_2$ and
$e_1+e_2+h_1+h_2$ belong to $D(A_0)$.
\end{remark}

\begin{proof}[Proof of Proposition~\refp{tfree}]
Fix two distinct elements $d_1,d_2\in D(A)$ and consider the subgroup
$H:=\<d_1,d_2\>$ generated by $d_1$ and $d_2$.

Suppose, to begin with, that the edges of $\Gam(A)$ corresponding to $d_1$
and $d_2$ are incident; that is, there are $a,b_1,b_2\in A$ such that
$d_1=a+b_1$ and $d_2=a+b_2$.
It is easy to see that the coset $a+H$ contains exactly three elements of $A$
(namely $b_1$, $b_2$ and $a$), while every other coset of $H$ contains at
most two elements of $A$---both conclusions in view of $d_1,\,d_2\in D(A)$.
Thus, the assumption $|A|\geq 2^{r-2}+3$ implies that there is a coset
containing exactly two elements of $A$. These two elements cannot differ by
$d_1$ or $d_2$ (again, since $d_1,\,d_2\in D(A)$); therefore they differ by
$d_1+d_2$, yielding a representation of $d_1+d_2$ as a sum of two elements of
$A$. Another representation is $d_1+d_2=b_1+b_2$, and the existence of two
representations shows that $d_1+d_2\notin D(A)$. The first assertion follows
since if $\Gam(A)$ were containing a triangle with two legs corresponding to
$d_1$ and $d_2$, then the third leg would correspond to $d_1+d_2$.

Assuming now that the edges of $\Gam(A)$ corresponding to $d_1$ and $d_2$ are
\emph{not} incident, find $a_1,a_2,b_1,b_2\in A$ such that $d_1=a_1+b_1$ and
$d_2=a_2+b_2$.
(Note that $a_1,\,a_2,\,b_1$ and $b_2$ are all distinct.) Then there are two
cosets of $H$ intersecting the set $\{a_1,a_2,b_1,b_2\}$. Each of these
cosets contains exactly two elements of $A$, while every other coset of $H$
contains at most two elements of $A$. If $|A|>2^{r-2}+3$, then there are at
least two cosets disjoint with $\{a_1,a_2,b_1,b_2\}$ and containing two
elements of $A$. This yields two distinct representations of $d_1+d_2$,
leading, as above, to the conclusion $d_1+d_2\notin D(A)$ and proving the
second assertion.
\end{proof}

Given a set $A\seq\Fr$, for each $a\in A$, we use $\deg(a)$ to denote the
degree of the vertex $a$ in $\Gam(A)$. Yet another fundamental property of
the unique representation graph is established by the following result.

\begin{proposition}\label{p:line}
Let $r\ge 1$ be an integer and suppose that $A\seq\Fr$ satisfies
$|A|>2^{r-2}+3$. If $(a_1,a_2)$ is an edge in $\Gam(A)$, then
  $$ \deg(a_1) + \deg(a_2) \ge |A| + |D(A)| - 2^{r-1}. $$
\end{proposition}

We present two different proofs.
\begin{proof}[First proof of Proposition~\refp{line}]
Let $A'$ denote the set of those elements of $A$ neighboring neither $a_1$
nor $a_2$ in $\Gam(A)$; thus, $|A'|=|A|-\deg(a_1)-\deg(a_2)$ by Proposition
\refp{tfree}. Then the sets
  $$ a_1+A',\ a_2+A',\ D(A),\ \text{and}\ a_1+a_2+D(A) $$
are easily seen to be pairwise disjoint, with the fact that the last two are
disjoint following from $a_1+a_2\in D(A)$ and Proposition \refp{tfree}, the
fact that the first two are disjoint following from $a_1+a_2\in D(A)$, and
the rest following from the definition of $A'$. Hence
  $$ 2^r \ge 2|A'| + 2|D(A)|
                      = 2 \big( |A| + |D(A)| - \deg(a_1)-\deg(a_2) \big). $$
\end{proof}

\begin{proof}[Second proof of Proposition~\refp{line}]
Since $a_1+a_2\in D(A)$ and the set $D(A)$ is sum-free by Proposition
\refp{tfree}, it contains at most one element from each coset of the
two-element subgroup $\<a_1+a_2\>$. On the other hand, $D(A)$ has exactly
$\deg(a_1)+\deg(a_2)-2$ elements in common with the set
$\{a_1,a_2\}+(A\stm\{a_1,a_2\})$, the size of which is $2(|A|-2)$, and which
is a union of cosets of $\<a_1+a_2\>$. It follows that
\begin{multline*}
 |D(A)|\le (\deg(a_1)+\deg(a_2)-2) + \frac12\,\big( 2^r-2(|A|-2) \big) \\
                                    = \deg(a_1) + \deg(a_2) + 2^{r-1} - |A|.
\end{multline*}
\end{proof}

We conclude this section deducing Theorem \reft{main} from Theorem
\reft{round}. To this end, we first derive from Proposition \refp{tfree} an
interesting property of sum-free sets. Thinking projectively, if $A$ is a
large cap in $\PG(r-1,2)$ and the point $p\notin A$ lies on the line
determined by a pair of points in $A$, then in fact there are \emph{many}
pairs of points in $A$ determining a line through $p$. (We remark that, for a
\emph{generic} subset of $\PG(r-1,2)$, not assumed to be a cap, the same
conclusion requires a much stronger assumption; cf.~Lemma~\refl{php-result}.)

In the definitions of a round set and the set $D(A)$ given above in this
section, we consider \emph{unordered} representations of elements of $\Fr$,
that is, representations which differ by the order of summands are considered
identical. This convention is extended onto the following corollary.

\begin{corollary}\label{c:sfnotround}
Let $r,\kap\ge 2$ be integers and suppose that $S\seq\Fr$ is a sum-free set
with $|S|>2^{r-2}+\kap$. Then every element of the sumset $2S$ has at least
$\kap$ representations (distinct under permutation of summands) as a sum of
two elements of $S$.
\end{corollary}

\begin{proof}
Assuming that an element $c\in 2S$ has fewer than $\kap$ representations as a
sum of two elements from $S$, we find a subset $S_0\seq S$ with $|S_0|\ge
|S|-(\kap-2)$ such that $c$ has exactly one representation as a sum of two
elements of $S_0$.

Let $A:=S_0\cup\{0\}$. Since $|S|>2^{r-2}+\kap$, we have $|A|>2^{r-2}+3$, so
in view of $S_0\seq D(A)$ and Proposition \refp{tfree}, we get $(2S_0)\cap
D(A)=\est$. Thus, every element of $2S_0$ has at least two representations as
a sum of two elements from $A$, and therefore at least two representations as
a sum of two elements from $S_0$ (since $(2S_0)\cap S_0=\est$), contradicting
the choice of $S_0$.
\end{proof}

\begin{proof}[Deduction of Theorem \reft{main} from Theorem \reft{round}]
As we have already observed, if $S\seq\Fr$ is a maximal sum-free set and
$s\in S\cup\{0\}$, then $\big(s+(S\cup\{0\})\big)\stm\{0\}$ is a minimal
$1$-saturating set. Suppose now that $r\ge 1$ is an integer and
$A\seq\Fr\stm\{0\}$ is a minimal $1$-saturating set with
$|A|>\frac{11}{36}\,2^r+3$. By Lemma \refl{1S2R}, either $A\cup\{0\}$ or $A$
is round. We show that, in the former case, $A$ is of the form required,
while the latter case cannot occur.

If $A\cup\{0\}$ is round, then by Theorem \reft{round} there exist a sum-free
set $S\seq\Fr$ and an element $g\in\Fr$ such that
$A\cup\{0\}=g+(S\cup\{0\})$. From $0\in g+(S\cup\{0\})$, it follows that
 $g\in S\cup\{0\}$, and $\Fr=2(A\cup\{0\})=2(S\cup\{0\})=S\cup 2S$ implies
that $S$ is a \emph{maximal} sum-free set (as remarked in Section 1), proving
the assertion in this case.

Suppose now that $A$ is round, so that by Theorem \reft{round} there exist a
sum-free set $S\seq\Fr$ and an element $g\in\Fr$ with $A=g+(S\cup\{0\})$. In
view of the previous paragraph, we may assume that
$A\cup\{0\}=g+(S\cup\{0,g\})$ is \emph{not} round, whence $S\cup\{g\}$ is not
sum-free (see the comment just above Theorem \ref{t:round}); that is, $g\in
2S$, and we write $g=s_1+s_2$ with $s_1,s_2\in S$. Notice that $0\notin A$
yields $g\ne 0$ and thus $s_1\ne s_2$, and that $2A=S\cup 2S$ and
  $$ A\cup 2A=S\cup(g+S)\cup 2S. $$

Let $S_1:=S\stm\{s_1\}$ and $A_1:=g+(S_1\cup\{0\})$. Since
  $$ |S| = |A|-1 > 2^{r-2}+2, $$
it follows from Corollary \refc{sfnotround} (applied with $\kap=2$) that
$2S_1=2S$. Consequently,
  $$ A_1\cup 2A_1=S_1\cup(g+S_1)\cup 2S. $$
On the other hand, as $g=s_1+s_2$ with $s_1\neq s_2$, we have
 $s_1,s_2\in S_1\cup(g+S_1)$, implying $S_1\cup(g+S_1)=S\cup(g+S)$; therefore,
$A_1\cup 2A_1=A\cup 2A$, contradicting the minimality of $A$.
\end{proof}

\section{Notation and Auxiliary Results.}\label{s:aux}

In this section, we deviate slightly from the flow of the proof to introduce
some important notation and results, preparing the ground for the rest of the
argument. We start with an easy consequence of the pigeonhole principle;
see, for instance, \cite[Lemma~2.1]{b:n} or \cite[Lemma~5.29]{b:gh}.
\begin{lemma}\label{l:php-result}
Let $B$ and $C$ be non-empty subsets of a finite abelian group $G$. If
$|B|+|C|\ge |G|+\kap$ with an integer $\kap\ge 1$, then every element of $G$
has at least $\kap$ representations as a sum of an element from
$B$ and an element from $C$.
\end{lemma}

We remark that, in Lemma \refl{php-result} and in the vast majority of
situations below, we consider representations of elements of $\Fr$ as sums of
elements from two \emph{potentially distinct} sets; therefore (in contrast
with Section~\refs{round}), representations are considered ordered.

Given a subgroup $H$ of an abelian group $G$, by $\phi_H$ we denote the
canonical homomorphism from $G$ onto the quotient group $G/H$.

For a subset $B$ of an abelian group $G$, the (maximal) period of $B$ will be
denoted by $\pi(B)$; recall that this is the subgroup of $G$ defined by
  $$ \pi(B) := \{g\in G\colon B+g=B \}, $$
and that $B$ is called \emph{periodic} if $\pi(B)\neq\{0\}$ and
\emph{aperiodic} otherwise. Thus, $B$ is a union of $\pi(B)$-cosets, and
$\pi(B)$ lies above every subgroup $H\le G$ such that $B$ is a union of
$H$-cosets. Observe also that $\pi(B)=G$ if and only if either $B=\est$ or
$B=G$, and that $\phi_{\pi(B)}(B)$ is an aperiodic subset of the group
$G/\pi(B)$.

\begin{theorem}[Kneser, \cite{b:kn1,b:kn2}; see also
                       \cite{b:m,b:n,b:gh}]\label{t:kneser}
Let $B$ and $C$ be finite, non-empty subsets of an abelian group $G$. If
  $$ |B+C| \le |B|+|C|-1, $$
then, letting $H:=\pi(B+C)$, we have
  $$ |B+C| = |B+H|+|C+H|-|H|. $$
\end{theorem}

\begin{corollary}\label{c:alldisjoint}
Let $r\ge 1$ be an integer and suppose that the sets $B,C\seq\Fr$ are
disjoint and non-empty. If $|B|+|C|>2^{r-1}$, then $B\cup C$ is not disjoint
with $B+C$.
\end{corollary}

\begin{remark}
If the elements $e_1,e_2\in\Fr$ and the subgroup $H<\Fr$ of index $4$ are so
chosen that $\Fr=\<e_1,e_2\>\oplus H$, then the sets $B:=e_1+H$ and
$C:=e_2+H$ are disjoint, and so are their union $B\cup C=\{e_1,e_2\}+H$ and
sumset $B+C=e_1+e_2+H$; at the same time, $|B|+|C|=2^{r-1}$. This shows that
the bound $2^{r-1}$ in Corollary \refc{alldisjoint} is sharp.
\end{remark}

\begin{proof}[Proof of Corollary \refc{alldisjoint}]
We proceed by induction on $r$. The case $r=1$ is immediate, and so we assume
 $r\ge 2$. Assuming, furthermore, that $B\cup C$ and $B+C$ are
disjoint, whereas $|B|+|C|>2^{r-1}$, we derive
  $$ |B+C|\le 2^r-|B|-|C| < |B|+|C|-1. $$
Set $H:=\pi(B+C)$. By Theorem \reft{kneser}, the subgroup $H$ is non-trivial
and
  $$ |(B+H)\stm B| + |(C+H)\stm C|
                             = |B+C| - |B| - |C| + |H| < |H| - 1. $$
The left-hand side can be interpreted as the total number of ``$H$-holes'' in
$B$ and $C$, showing that $B+H$ and $C+H$ are disjoint (since $B$ and $C$ are
themselves disjoint). By the same reasoning, these two sets are also disjoint
with $B+C$ (as $B+C$ is disjoint with both $B$ and $C$, and $\pi(B+C)=H$, so
there are no ``$H$-holes'' in $B+C$). Consequently, $\phi_H(B)$ and
$\phi_H(C)$ are disjoint, non-empty subsets of the group $\Fr/H$, and
$\phi_H(B)\cup\phi_H(C)$ is disjoint with $\phi_H(B)+\phi_H(C)=\phi_H(B+C)$.
This contradicts the induction hypothesis in view of
\begin{multline*}
  |\phi_H(B)|+|\phi_H(C)| = (|B+H|+|C+H|)/|H| \\
                         \ge (|B|+|C|)/|H| > 2^{r-1}/|H| = \frac12\,|\Fr/H|.
\end{multline*}
\end{proof}

For an integer $k$ and subsets $B$ and $C$ of an additively written group,
let $B\multplus[k]C$ denote the set of all those group elements with at least
$k$ representations as $b+c$ with $b\in B$ and $c\in C$; thus, for instance,
$B\multplus[1]C=B+C$. We need a corollary of the following theorem, which is
(a refinement of) a particular case of the main result of \refb{g}.

\begin{theorem}[Grynkiewicz, \protect{\cite[Theorem 1.2]{b:g}}]\label{t:g}
Let $G$ be an abelian group and suppose that $B,C\seq G$ are finite and
satisfy $\min\{|B|,|C|\}\ge 2$. Then either
  $$ |B\multplus[1]C| + |B\multplus C| \ge 2|B| + 2|C| - 4, $$
or there exist subsets $B'\seq B$ and $C'\seq C$ with
\begin{gather*}
  l := |B\stm B'| + |C\stm C'| \le 1, \\
  B'+C' = B'\multplus C' = B\multplus C,
\end{gather*}
and
\begin{align*}
  |B\multplus[1]C| + |B\multplus C|
    &\ge 2|B| + 2|C| - (2-l)(|H|-\rho) - 2l \\
    &\ge 2|B|+2|C|-2|H|,
\end{align*}
where $H=\pi(B\multplus C)$ and $\rho=|(B'+H)\stm B'|+|(C'+H)\stm C'|$.
\end{theorem}
(For our present purposes, the reader can completely ignore the definitions
of $H$ and $\rho$ in the statement of Theorem \ref{t:g} and the part of the
conclusion involving these quantities.)

\begin{corollary}\label{c:S2}
If $G$ is a finite abelian group and $B,C\seq G$ satisfy
 $\min\{|B|,|C|\}\ge 2$, then
  $$ |B\multplus C| \ge \min \{ 2|B|+2|C|-4-|G|, |B|-1 \}. $$
\end{corollary}

\begin{proof}
If $|B\multplus[1]C|+|B\multplus C|\ge 2|B|+2|C|-4$, then
 $|B\multplus C|\ge 2|B|+2|C|-4-|G|$ follows trivially. Otherwise, we apply
Theorem~\reft{g} to find $B'\seq B$ and $C'\seq C$ satisfying
 $|B\stm B'|+|C\stm C'|\le 1$ and $B'+C'=B'\multplus C'=B\multplus C$. Now
  $$ |B\multplus C| = |B'+C'| \ge |B'| \ge |B|-1. $$
\end{proof}

Finally, we prove several simple graph-theoretic lemmas and apply them to the
unique representation graph.

Recall that the matching number of a graph is the largest number of edges in
a matching of the graph.

\begin{lemma}\label{l:matching2}
Let $(V,E)$ be a triangle-free graph without isolated vertices, such that the
matching number of $(V,E)$ does not exceed $2$. If $|V|\ge 6$, then $(V,E)$
is either a star or a union of two stars, possibly with an edge between their
centers. More precisely, there is a partition
 $V=\{v_1,v_2\}\cup V_0\cup V_1\cup V_2$ such that $E$ consists of all pairs
$(v_1,v)$ with $v\in V_0\cup V_1$, all pairs $(v_2,v)$ with
 $v\in V_0\cup V_2$, and, possibly, the pair $(v_1,v_2)$.
\end{lemma}

\begin{proof}
We notice that $(V,E)$ does not contain a pentagon: for otherwise, one could
construct a matching of size $3$ using two edges of the pentagon and an edge
incident with a vertex outside the pentagon. Furthermore, $(V,E)$ does not
contain cycles of length $6$ or more. Consequently, $(V,E)$ contains no odd
cycles; hence it is bipartite.

As a result, by K\"onig's theorem, $(V,E)$ has a vertex cover of size at most
$2$. Now if $\{v\}$ is a vertex cover, then $(V,E)$ is a star with the center
at $v$, and if $\{v_1,v_2\}$ with $v_1\neq v_2$ is a vertex cover, then the
assertion follows by letting $V_0$ be the set of common neighbors of $v_1$
and $v_2$, and, for $i\in\{1,2\}$, defining $V_i$ to be the set of all
neighbors of $v_i$ in $V\stm (V_0\cup\{v_1,v_2\})$.
\end{proof}

\begin{lemma}\label{l:linegraph}
Let $\del\ge 1$ be an integer and suppose that $(V,E)$ is a graph such that
$\deg(v_1)+\deg(v_2)\ge\del$ holds for every edge $(v_1,v_2)\in E$. If
$(V,E)$ has no isolated vertices, then $|E|\ge(1-\del^{-1})|V|$.
\end{lemma}

\begin{remark}
Equality is attained if $(V,E)$ is a disjoint union of stars with $\delta$
vertices each.
\end{remark}

\begin{proof}[Proof of Lemma \refl{linegraph}]
For $\del\le 2$, the assertion is immediate. Assume therefore that $\del\ge
3$ and, for each $j\in[1,\del-2]$, let $V_j:=\{v\in V\colon\deg(v)=j\}$;
also, let $V_+:=\{v\in V\colon\deg(v)\ge \del-1\}$, so that $V$ is the
disjoint union of $V_1\longc V_{\del-2}$ and $V_+$. Evidently, we have
\begin{align}
  \sum_{v\in V_j} \deg(v) &= j|V_j|,\quad j\in[1,\del-2], \label{e:linegraph1}
\intertext{and}
  \sum_{v\in V_+} \deg(v) &\ge (\del-1)|V_+|. \label{e:linegraph2}
\intertext{Also,}
  \sum_{v\in V_+} \deg(v) &\ge |V_1|,\label{e:linegraph3}
\end{align}
as every vertex from $V_1$ is adjacent to a vertex from $V_+$ in view of the
hypothesis $\deg(v_1)+\deg(v_2)\ge\del$. Taking the sum of inequality
\refe{linegraph2} with weight $2\del^{-1}$, inequality \refe{linegraph3} with
weight $1-2\del^{-1}$, and equations \refe{linegraph1} with weight $1$ for
each $j\in[1,\del-2]$, we get
  $$ 2|E| = \sum_{v\in V} \deg(v)
       \ge (2-2\del^{-1})|V_1| + \sum_{j=2}^{\del-2} j|V_j|
                                                        + (2-2\del^{-1})|V_+|
         \ge (2-2\del^{-1})|V|. $$
\end{proof}

Applying Lemma~\refl{linegraph} with $\del=2$ to the unique representation
graph of a round set, we get the following corollary.
\begin{corollary}\label{c:DishalfA}
If $r\ge 1$ is an integer and $A\seq\Fr$ is a round set, then
$|D(A)|\ge\frac12|A|$.
\end{corollary}

\begin{lemma}\label{lemma-highfive}
Let $t$ be the matching number of a graph $(V,E)$. If $(V,E)$ does not have
isolated vertices, then
  $$ |V| \le |E|+t. $$
\end{lemma}

\begin{proof}
If $T$ is a matching with $|T|=t$ edges, then $(V,E)$ has $|V|-2t$ vertices
not incident with the edges of $T$. By the maximality of $T$, no two of these
vertices are adjacent, and thus each of them is incident to an edge from
$E\stm T$, since $(V,E)$ contains no isolated vertices. Consequently,
  $$ |E|  \ge t + (|V|-2t) = |V| - t. $$
\end{proof}

Since the matching number of a graph does not exceed the number of edges in
the graph, the following corollary strengthens Corollary \refc{DishalfA}.
\begin{corollary}\label{c:matching}
If $r\ge 1$ is an integer and $A\seq\Fr$ is a round set, then
$|A|\le|D(A)|+t$, where $t$ is the matching number of $\Gam(A)$.
\end{corollary}

\section{A ``Light Version'' of Theorem \reft{main}.}\label{s:light}

In this section, we combine the tools, developed so far, to prove the
following, slightly weaker version of Theorem \reft{main}.
\renewcommand{\theprimetheorem}{\reft{main}$'$}
\begin{primetheorem}\label{t:mainlight}
Let $r\ge 1$ be an integer. A set $A\seq\Fr\stm\{0\}$ with
$|A|>\frac13\,2^r+2$ is minimal $1$-saturating if and only if there are a
maximal sum-free set $S\seq\Fr$ and an element $s\in S\cup\{0\}$ such that
$A=\big(s+(S\cup\{0\})\big)\stm\{0\}$.
\end{primetheorem}

Examining the deduction of Theorem \reft{main} from Theorem~\reft{round} at
the end of Section \refs{round}, the reader will see that, in an identical
way, Theorem \reft{mainlight} can be obtained
from the following ``week version'' of Theorem \reft{round} (the only
difference being that for Theorem \reft{main}$'$ one needs the estimate
$|A|\le 2^{r-1}$, mentioned in Section 1, to exclude the cases $r\le 3$ where
$|A|>\frac13\,2^r+2$ does not imply $|A|>2^{r-2}+3$).
\renewcommand{\theprimetheorem}{\reft{round}$'$}
\begin{primetheorem}\label{t:roundlight}
Let $r\ge 1$ be an integer and suppose that $A\seq\Fr$ is round. If
$|A|>\frac13\,2^r+2$, then there is a sum-free set $S\seq\Fr$ and an element
$g\in\Fr$ such that $A=g+(S\cup\{0\})$.
\end{primetheorem}

Thus, all we need is to prove Theorem~\reft{roundlight}.

\begin{proof}[Proof of Theorem~\reft{roundlight}]
Suppose that $|A|>\frac13\,2^r+2$. As mentioned in Section~\refs{intro}, the
size of a round set in $\Fr$ does not exceed $2^{r-1}$. Consequently, the
hypotheses imply $r\ge 4$ and, furthermore, $|A|>2^{r-2}+3$; this is
implicitly used below to invoke Propositions~\refp{tfree} and~\refp{line}.

Set $\del:=|A|+|D(A)|-2^{r-1}$. By Corollary~\refc{DishalfA}, we have
  $$ \del \ge \frac32\,|A| - 2^{r-1} > 0; $$
thus, Proposition~\refp{line} and Lemma~\refl{linegraph} give
$\del|D(A)|\ge(\del-1)|A|$. Substituting the value of $\del$ and rearranging
the terms, we rewrite this estimate as
\begin{equation}\label{e:what}
  f(|D(A)|) \le |A|(2^{r-1}+1-|A|),
\end{equation}
where $f$ is the real function defined by $f(x):=x(2^{r-1}-x)$.

Since $f$ is concave, $|D(A)|\ge\frac12\,|A|$ by Corollary~\refc{DishalfA},
and
  $$ \min \{ f(|A|/2),\,f(|A|-2) \} > |A|(2^{r-1}+1-|A|) $$
(which follows by a straightforward computation using the assumption on the
size of $A$), we derive from \refe{what} that
\begin{equation}\label{e:lightDlarge}
  |D(A)|\ge|A|-1.
\end{equation}

In view of Lemma~\refl{sf=star}, it suffices to show that $\Gam(A)$ has a
spanning star; that is (since $\Gam(A)$ is triangle-free by
Proposition~\refp{tfree} and has no isolated vertices), that the matching
number of $\Gam(A)$ is equal to $1$. Suppose, for a contradiction, that
$\Gam(A)$ has a two-edge matching $T$. By Proposition~\refp{line}, incident
with each of the two edges of $T$ are $\del-2$ edges of $\Gam(A)$. Moreover,
since $\Gam(A)$ is triangle-free, there are at most two edges of $\Gam(A)$
incident with both edges of $T$. Consequently, the total number of edges of
$\Gam(A)$ is at least
  $$ |D(A)| \ge 2(\del-2) + |T| - 2 = 2|A| + 2|D(A)| - 2^r - 4. $$
Using \refe{lightDlarge} and the assumption $|A|>\frac13\,2^r+2$, we derive
  $$ 2^r + 4 \ge 2|A| + |D(A)| \ge 3|A| - 1 > 2^r + 5, $$
a contradiction.
\end{proof}

\section{Securing Two Isolated Edges.}\label{s:twoisolated}

In this section, we prove Theorem \reft{round} under the extra assumption
that $\Gam(A)$ has at most one isolated edge; the case where $\Gam(A)$ has
two or more isolated edges is dealt with in Sections \refs{2edgesbasics} and
\refs{2edgescompletion}. We split the argument into two lemmas.

\begin{lemma}\label{l:twoisolated}
Let $r\ge 1$ be an integer and suppose that $A\seq\Fr$ is round. If $\Gam(A)$
has at most one isolated edge and $|A|>\frac{3}{5}\, 2^{r-1}+\frac{13}{5}$,
then the matching number of $\Gam(A)$ is at most $2$.
\end{lemma}

The proof is a minor modification of that of Theorem \reft{roundlight}.

\begin{proof}[Proof of Lemma~\refl{twoisolated}]
Since $|A|\le 2^{r-1}$, we have $r\ge 4$, and thus
$\frac{3}{5}\,2^{r-1}+\frac{13}{5}>2^{r-2}+3$; this will allow us to apply
Propositions~\refp{tfree} and~\refp{line}.

If $\Gam(A)$ does not have isolated edges then, applying Lemma
\refl{linegraph} to the graph $\Gam(A)$, we get $|D(A)|\ge \frac23\,|A|$; if
$\Gam(A)$ has one isolated edge, then, applying Lemma \refl{linegraph} to the
graph $\Gam(A)$ with this edge removed, we get $|D(A)|\ge
1+\frac23\,(|A|-2)=\frac23\,|A|-\frac13$. In any case, letting
$\del:=|A|+|D(A)|-2^{r-1}$ and assuming
$|A|>\frac{3}{5}\,2^{r-1}+\frac{13}{5}$, we have
  $$ \del \ge \frac53\,|A| - 2^{r-1} - \frac13 > 0. $$
Consequently, applying Proposition \refp{line} and Lemma \refl{linegraph}, we
obtain $\del|D(A)|\ge (\del-1)|A|$. Substituting the value of $\del$,
rearranging the terms, and letting $f(x):=(2^{r-1}-x)x$, we re-write this
estimate as
  $$ f(|D(A)|) \le (2^{r-1}+1-|A|)|A|. $$

We notice that $f(x)$ is concave, that
\begin{align*}
  f\(\frac23\,|A|-\frac13\)
    &= \frac23\,|A| \( 2^{r-1} - \frac23\,|A| + \frac13 \)
                         - \frac13\, \( 2^{r-1} - \frac23\,|A| + \frac13 \) \\
    &> \( \frac23\,2^{r-1} - \frac49\,|A| + \frac29 \) |A| - \frac13\,|A| \\
    &= \( \frac23\,2^{r-1} - \frac49\,|A| - \frac19 \) |A| \\
    &> (2^{r-1}+1-|A|)|A|,
\intertext{and that}
  f(|A|-3)
    &= (2^{r-1}-|A|+3)(|A|-3) \\
    &= (2^{r-1}-|A|+1)|A| + 5|A| - 3\cdot 2^{r-1} - 9 \\
    &> (2^{r-1}-|A|+1)|A|,
\end{align*}
where all three estimates follow from
$|A|>\frac{3}{5}\,2^{r-1}+\frac{13}{5}$. Thus, in view of
$|D(A)|\ge\frac23\,|A|-\frac13$, we conclude that, indeed,
\begin{equation}\label{e:duckhat}
  |D(A)| \ge |A|-2.
\end{equation}

Suppose now by contradiction that $\Gam(A)$ possesses a three-edge matching
$T$. Using Proposition \refp{line} to count the edges of $\Gam(A)$ incident
to those in $T$, and also taking into account the three edges of $T$, we get
  $$ |D(A)| \ge 3 \big(|A|+|D(A)|-2^{r-1}-2\big) + 3 - 6; $$
for any edge incident to two different edges from $T$ joins two vertices from
$T$ while, since $\Gam(A)$ is triangle-free (by Proposition \refp{tfree}) and
$T$ is a matching in $\Gam(A)$, the graph induced by the six vertices of $T$
has at most six edges not in $T$. Rearranging the terms and applying
\refe{duckhat} gives
  $$ 2^{r-1} \ge |A| + \frac23\,|D(A)| - 3 \ge \frac53\,|A| - \frac{13}{3}, $$
which contradicts the assumption on $|A|$.
\end{proof}

\begin{lemma}\label{l:twostars}
Let $r\ge 1$ be an integer. If $A\seq\Fr$ is a round set with
$|A|>2^{r-2}+3$, then the matching number of $\Gam(A)$ is distinct from $2$.
\end{lemma}

\begin{proof}
Assume for a contradiction that $A\seq\Fr$ is round, $|A|>2^{r-2}+3$, and the
matching number of $\Gam(A)$ is equal to $2$.
From $2^{r-1}\ge|A|>2^{r-2}+3$ we derive $2^{r-2}>3$, and then $|A|>6$.
Hence, by Proposition \refp{tfree} and Lemma \refl{matching2}, there exist
distinct elements $a_1,a_2\in A$ and disjoint subsets $A_0,A_1,A_2\seq
A\stm\{a_1,a_2\}$ such that
 $A=\{a_1,a_2\}\cup A_0\cup A_1\cup A_2$ and
  $$ D(A)\stm\{a_1+a_2\} = \big( a_1+(A_1\cup A_0) \big)
                                        \cup \big( a_2+(A_2\cup A_0) \big). $$
Indeed, $a_1+a_2\in D(A)$ holds: else, for some
 $a',a''\in A_0\cup A_1\cup A_2$, we would have $a_1+a_2=a'+a''$,
contradicting the fact that either $a_1+a'$ or $a_2+a'$ is uniquely
representable (up to permutation of summands) as a sum of two elements of
$A$.

By Proposition \refp{tfree}, $\Gam(A)$ is triangle-free, and consequently,
$A_0=\est$: for $a_1$ and $a_2$ are joined by an edge in $\Gam(A)$ and
therefore have no common neighbors. Hence,
  $$ D(A) = \{a_1+a_2\} \cup (a_1+A_1) \cup(a_2+A_2), $$
where the union is disjoint by the definition of $D(A)$. For $i\in\{1,2\}$,
we write $D_i:=a_i+A_i$, and we consider the sets $B_i:=\{0,a_1+a_2\}+D_i$.
Since $a_1+a_2\in D(A)$ and $D_1,D_2\seq D(A)$, and since $D_1\cap D_2=\est$,
it follows in view of Proposition \refp{tfree} that $B_1\cap B_2=\est$ and
  $$ |B_1|+|B_2| = 2|D_1|+2|D_2| = 2(|A|-2) > 2^{r-1}. $$

We claim now that the sumset $B_1+B_2=\{0,a_1+a_2\}+D_1+D_2$ is disjoint from
the union $B_1\cup B_2=\{0,a_1+a_2\}+(D_1\cup D_2)$; which, since
$\{0,a_1+a_2\}$ is a subgroup, is equivalent to $\{0,a_1+a_2\}+D_1+D_2$ being
disjoint with $D_1\cup D_2$. To see this, assume that $\{0,a_1+a_2\}+D_1+D_2$
is not disjoint with, say, $D_1$. As $(D_1+D_2)\cap D_1=\est$ by Proposition
\refp{tfree}, this assumption yields $(a_1+a_2+D_1+D_2)\cap D_1\neq\est$;
that is, $a_1+a_2+d_1+d_2=d_1'$ for some $d_1,d_1'\in D_1$ and $d_2\in D_2$.
Letting $\alp_i:=a_i+d_i\ (i\in\{1,2\})$, we re-write this equality as
$\alp_1+\alp_2=d_1'$ and obtain a contradiction observing that $\alp_1\in
A_1\seq A\stm\{a_1\}$ and $\alp_2\in A_2\seq A\stm\{a_1\}$, whereas $d_1'\in
D_1$ shows that the only representation of $d_1'$ as a sum of two elements of
$A$ involves $a_1$ as a summand.

Applying Corollary \refc{alldisjoint} to the sets $B_1$ and $B_2$, we
conclude that one of them is empty. Consequently, either $A_1$ or $A_2$ is
empty. Thus, $\Gam(A)$ is a star, whence the matching number of $\Gam(A)$ is
$1$, contrary to an assumption at the beginning of the proof.
\end{proof}

\section{Using Two Isolated Edges: The Coset Structure.}
  \label{s:2edgesbasics}

As follows from Lemmas \refl{sf=star}, \refl{twoisolated} and
\refl{twostars}, and since $\frac{11}{36}>\frac{3}{10}$, to complete the
proof of Theorem~\reft{round}, it remains to consider the case where
$\Gam(A)$ has at least two isolated edges. Accordingly, we assume in this and
the next section that $r\ge 1$ is an integer and that $A\seq\Fr$ is a round
set such that $\Gam(A)$ has two (or more) isolated edges, and we show that
$|A|<\frac{11}{36}\,2^r+3$.

Shifting $A$, if necessary, we assume that $0\in A$ and that $a_1,\,a_2$, and
$a_3$ are elements of $A$, distinct from $0$ and each other, such that
$(0,a_1)$ and $(a_2,a_3)$ are isolated edges of $\Gam(A)$. We consider the
subgroups $L=\<a_1,a_2,a_3\>,\ K^-=\<a_3,a_1+a_2\>,\ K^+=\<a_2,a_1+a_3\>$,
and $H=\<a_1+a_2+a_3\>$; thus,
  $$ |L|=8,\ |K^-|=|K^+|=4,\ H = K^-\cap K^+,\ \text{and}\ |H|=2. $$
Our argument is based on a careful study of the distribution of the elements
of $A$ and $D(A)$ in the cosets of $L$. The goal of the present section is to
establish some basic facts about this distribution.

For $g\in\Fr$, we write $A_g:=A\cap(g+L)$ and $D_g:=D(A)\cap(g+L)$.
Evidently, we have $\{0,a_1,a_2,a_3\}\seq A_0$, and it is easy to see
that, indeed, $A_0=\{0,a_1,a_2,a_3\}$. Next,
from $\{a_1,a_2+a_3\}\seq(2A_0)\cap D(A)$, it follows that
\begin{equation}\label{e:coset2A}
  (2A_g)\cap \{a_1,a_2+a_3\} = \est
\end{equation}
for $g\notin L$, and the fact that $(0,a_1)$ and $(a_2,a_3)$ are isolated
edges gives
\begin{equation}\label{e:cosetAD}
  (A_g+A_0) \cap D_g = \est,
\end{equation}
under the same assumption.
Furthermore, in view of Proposition \refp{tfree} and since
$a_1,\,a_2+a_3\in D(A)$, we have
\begin{equation}\label{e:coset2D}
  (2D_g) \cap \{a_1,a_2+a_3\} = \est,
\end{equation}
for each $g\in\Fr$.

An immediate corollary of \refe{coset2A} and Lemma \refl{php-result} is that
$|A_g|\le 4$ holds for every element $g\in\Fr$. With this in mind, for
$g\in\Fr$ and $i\in[0,4]$, we say that the coset $g+L$ is of type $i$ if
$|A_g|=i$, and we denote by $n_i$ the number of \emph{non-zero} $L$-cosets of
type $i$ (so that $L$ is not counted in $n_4$); hence,
\begin{equation}\label{e:sum-n_i-is-m/8}
  n_0+n_1+n_2+n_3+n_4 = 2^{r-3}-1
\end{equation}
and
\begin{equation}\label{e:sum-n_i-is-A}
  n_1 + 2n_2 + 3n_3 + 4n_4 = |A|-4.
\end{equation}

We now introduce a manner of pictorially representing the distribution of
subsets of $\Fr$ in $L$-cosets that will help elucidate the otherwise tedious
arguments needed for this section, and which may be helpful to keep in mind
for the next section as well. Specifically, given a set $X\seq\Fr$ and an
element $g\in\Fr$, we represent the elements of $X$ in the coset $g+L$ by a
diagram like
  $$ X\cap(g+L):
       \begin{array}{rc}
         \begin{array}{r}
           \vspace{.1cm}
           \begin{array}{r}
             \vspace{.1cm}
             \begin{array}{r}
               g \hspace{1.4cm} \longrightarrow \\
               g+a_1+a_2+a_3    \longrightarrow
             \end{array} \\
             \vspace{.1cm}
             \begin{array}{r}
               g+a_1     \hspace{.9cm} \longrightarrow \\
               g+a_2+a_3 \hspace{.4cm} \longrightarrow
             \end{array}
           \end{array} \\
           \vspace{.1cm}
           \begin{array}{r}
             \vspace{.1cm}
             \begin{array}{r}
               g+a_2     \hspace{.9cm} \longrightarrow \\
               g+a_1+a_3 \hspace{.4cm} \longrightarrow
             \end{array} \\
             \vspace{.1cm}
             \begin{array}{r}
               g+a_1+a_2 \hspace{.4cm} \longrightarrow \\
               g+a_3     \hspace{.9cm} \longrightarrow
             \end{array}
           \end{array}
         \end{array}
         \hspace{-1cm} &
         \begin{array}{l}
           \vspace{.1cm}
           \left(
           \begin{array}{l}
             \vspace{.1cm}
             \left(
             \begin{array}{c}
               \bullet \\
               \bullet
             \end{array}
             \right) \\
             \vspace{.1cm}
             \left(
             \begin{array}{c}
               \circ \\
               \circ
             \end{array}
             \right)
           \end{array}
           \right) \\
           \vspace{.1cm}
           \left(
           \begin{array}{c}
             \vspace{.1cm}
             \left(
             \begin{array}{l}
               \circ \\
               \bullet
             \end{array}
             \right) \\
             \vspace{.1cm}
             \left(
             \begin{array}{c}
               \circ \\
               \circ
             \end{array}
             \right)
           \end{array}
           \right)
         \end{array}
       \end{array} $$
where each filled dot represents an element (as labeled) contained in $X$,
and each open dot represents an element not in $X$. Note that this
representation depends, though only up to translation, on the choice of the
element $g$ within the $L$-coset.

We remark that two blocks of points enclosed by parentheses of the same level
are cosets of the same subgroup; say, the four two-point blocks correspond to
the four $H$-cosets contained in $g+L$.

As an example, the distribution of the elements of $A$ in $L$ can be depicted
as
  $$
     A\cap L:
     \begin{array}{rc}
    \begin{array}{r}
     \vspace{.1cm}
  \begin{array}{r}
     \vspace{.1cm}
  \begin{array}{r}
    0 \hspace{1cm}\longrightarrow\\
    a_1+a_2+a_3 \longrightarrow
  \end{array}
 \\
    \vspace{.1cm}
  \begin{array}{r}
    a_1 \hspace{.9cm}\longrightarrow\\
    a_2+a_3\hspace{.4cm} \longrightarrow
  \end{array}
  \end{array}
 \\
     \vspace{.1cm}
  \begin{array}{r}
     \vspace{.1cm}
  \begin{array}{r}
    a_2\hspace{.9cm}\longrightarrow \\
    a_1+a_3\hspace{.4cm}\longrightarrow
  \end{array}
 \\
     \vspace{.1cm}
  \begin{array}{r}
    a_1+a_2\hspace{.4cm} \longrightarrow\\
    a_3\hspace{.9cm}\longrightarrow
  \end{array}
  \end{array}
  \end{array}
\hspace{-1cm}
    &
    \left. \begin{array}{l}
     \vspace{.1cm}
    \left. \(
  \begin{array}{l}
     \vspace{.1cm}
  \left. \(
\begin{array}{c}
    \bullet \\
    \circ
  \end{array}
\)\rfr H
 \\
    \vspace{.1cm}
    \(
  \begin{array}{c}
    \bullet \\
    \circ
  \end{array}
\)
  \end{array}
\) \rfr M \\
     \vspace{.1cm}
    \(
  \begin{array}{c}
     \vspace{.1cm}
    \(
  \begin{array}{l}
    \bullet \\
    \circ
  \end{array}
\)
\\
     \vspace{.1cm}
    \(
  \begin{array}{c}
    \circ \\
    \bullet
  \end{array}
\)
  \end{array}
\hspace{.8cm} \)
  \end{array} \rfr L
    \end{array}
  $$
where we have used braces to label the subgroups $H$, $L$ and
$M:=\<a_1,a_2+a_3\>$. Furthermore, $K^-$ and $K^+$ are located in $L$ as
follows:
 $$ \begin{array}{cc}
  K^-: \begin{array}{c}
     \vspace{.1cm}
    \left(
  \begin{array}{c}
     \vspace{.1cm}
    \left(
  \begin{array}{c}
    \bullet \\
    \bullet \\
  \end{array}
\right) \\
    \vspace{.1cm}
    \left(
  \begin{array}{c}
    \circ \\
    \circ \\
  \end{array}
\right) \\
  \end{array}
\right) \\
     \vspace{.1cm}
    \left(
  \begin{array}{c}
     \vspace{.1cm}
    \left(
  \begin{array}{c}
    \circ \\
    \circ \\
  \end{array}
\right) \\
     \vspace{.1cm}
    \left(
  \begin{array}{c}
    \bullet \\
    \bullet \\
  \end{array}
\right) \\
  \end{array}
\right) \\
  \end{array}
  &
  \hspace{1cm}
K^+: \begin{array}{c}
     \vspace{.1cm}
    \left(
  \begin{array}{c}
     \vspace{.1cm}
    \left(
  \begin{array}{c}
    \bullet \\
    \bullet \\
  \end{array}
\right) \\
    \vspace{.1cm}
    \left(
  \begin{array}{c}
    \circ \\
    \circ \\
  \end{array}
\right) \\
  \end{array}
\right) \\
     \vspace{.1cm}
    \left(
  \begin{array}{c}
     \vspace{.1cm}
    \left(
  \begin{array}{c}
    \bullet \\
    \bullet \\
  \end{array}
\right) \\
     \vspace{.1cm}
    \left(
  \begin{array}{c}
    \circ \\
    \circ \\
  \end{array}
\right) \\
  \end{array}
\right) \\
  \end{array}
\end{array}$$

With the above diagrams in mind, we see that (\ref{e:coset2A}) is just the
statement that any two elements of $A$ from the same $M$-coset, excepting the
two $M$-cosets contained in $L$, are actually from the same $H$-coset. Thus
$|(g+M)\cap A|\le 2$ for each $g\in\F_2^r$, and consequently, given any $g\in
\F_2^r\stm L$, we can find $x,y\in g+L$ (one element for each of the two
$M$-cosets contained in $g+L$) such that $A_g\seq(x+H)\cup(y+H)$. Since
$(x+H)\cup(y+H)$ is either a $K^+$ or $K^-$-coset for any choice of $x$ and
$y$ in the same $L$-coset, we conclude that $A_g$ is contained either in a
single $K^+$-coset, or in a single $K^-$-coset. Using the homomorphism
notation from Section~\refs{aux}, we record this observation as follows.

\begin{claim}\label{m:onecoset}
For every $g\in\Fr\stm L$, we have
$\min\{|\phi_{K^-}(A_g)|,|\phi_{K^+}(A_g)|\}\le 1$.
\end{claim}

Refining our classification of cosets of $L$, for $i\in[2,4]$ and $g\in\Fr$,
we say that the coset $g+L$ is of type $i^0$ if it is of type $i$ and, in
addition,
  $$ |\phi_{K^-}(A_g)|=|\phi_{K^+}(A_g)|=1; $$
that $g+L$ is of type $i^-$ if it is of type $i$ and, in addition,
  $$ |\phi_{K^+}(A_g)|>|\phi_{K^-}(A_g)|=1; $$
and finally, that $g+L$ is of type $i^+$ if it is of type $i$ and
  $$ |\phi_{K^-}(A_g)|>|\phi_{K^+}(A_g)|=1. $$
Let $n_i^0,n_i^-$, and $n_i^+$ denote the number of non-zero cosets of the
corresponding types. From this definition, Claim~\refm{onecoset}, and the
observation that if $|\phi_{K^-}(A_g)|=|\phi_{K^+}(A_g)|=1$, then
$|\phi_H(A_g)|=1$ and thus $|A_g|\le 2$, it follows that
  $$ n_2 = n_2^0+n_2^-+n_2^+,\ n_3^0 = n_4^0=0,\ n_3 = n^-_3+n_3^+,
                                          \ \text{and}\ n_4 = n_4^-+n_4^+. $$

\begin{claim}\label{m:Dg}
For every $g\in\Fr$, we have
\begin{align*}
  |D_g| &= \begin{cases}
               0,\ &\text{if $g+L$ is of type $2^0,\,3$, or $4$,
                                                  and $g\notin L$}; \\
               2,\ &\text{if $g\in L$};
             \end{cases}
\intertext{furthermore,}
  |D_g| &\le \begin{cases}
                 2,\ &\text{if $g+L$ is of type $1,2^-$, or $2^+$}; \\
                 4\ &\text{if $g+L$ is of type $0$}.
               \end{cases}
\end{align*}
\end{claim}

\begin{proof}
If $g+L$ is of type $2^0$, then $A_g$ is an $H$-coset as $|A_g|=2$ and
$A_g$ is contained in the intersection of a $K^-$-coset and a
$K^+$-coset. If $g\notin L$ and $g+L$ is of type $3$ or $4$, then $A_g$
has two elements in the same $M$-coset, hence $A_g$ contains an $H$-coset
by the above observation that two elements of $A_g$, falling into the
same $M$-coset, are actually in the same $H$-coset. As a result, if
$g\notin L$ and $g+L$ is of type $2^0,3$, or $4$, then $A_g$ contains an
$H$-coset, and without loss of generality we assume $g+H\seq A_g$.
However, $g+H+A_0=g+L$ (as is readily apparent from the diagram for
$A_0$), whence $A_0+A_g=g+L$ and thus \refe{cosetAD} implies $|D_g|=0$.

By \refe{coset2D} and since $\{a_1,a_2+a_3\}\seq D(A)$, the set $D(A)$ is
disjoint with $\{0,a_1+a_2+a_3\}$, and the assumption that the edge
$(a_2,a_3)$ is isolated shows that $D(A)$ is also disjoint with
$\{a_2,a_3,a_1+a_2,a_1+a_3\}$. (If, for instance, we had $a_2\in D(A)$,
then $a_2$ would be adjacent to $0$; if we had $a_1+a_2\in D(A)$, then
$a_2$ would be adjacent to $a_1$ etc.) Consequently, if $g\in L$, then
$D_g=\{a_1,a_2+a_3\}$ and thus $|D_g|=2$.

Next, if $g\in A$ and $g\notin L$, then by \refe{cosetAD} the set $D_g$
is disjoint with $A_0+A_g\supseteq\{g,g+a_1,g+a_2,g+a_2,g+a_3\}$. Also,
\refe{coset2D} shows that $D_g$ can possibly contain at most one of
$g+a_2+a_3$ and $g+a_1+a_2+a_3$, and similarly $D_g$ can possibly contain
at most one of $g+a_1+a_3$ and $g+a_1+a_2$. It follows that $|D_g|\le 2$
whenever $g\notin L$ and $g+L$ is not of type $0$.

Finally, the fact that $|D_g|\le 4$ for each $g\in\Fr$ is a direct
consequence of \refe{coset2D} and Lemma~\refl{php-result}.
\end{proof}

\begin{claim}\label{m:DinXY}
For every $g\in\Fr$ such that $g+L$ is of type $1$, $2^+$  or $2^-$, there
exists a subset $\tD_g\seq g+L$ with  $D_g\seq\tD_g$ and
$|\tD_g|=|\phi_H(\tD_g)|$; moreover,
\begin{itemize}
\item[(i)]   if $g+L$ is of type $1$, then
    $|\tD_g|=|\varphi_H(\tD_g)|=4$;
\item[(ii)]  if $g+L$ is of type $2^-$, then
    $|\tD_g|=|\varphi_H(\tD_g)|=2$ and $|\phi_{K^-}(\tD_g)|=1$;
\item[(iii)] if $g+L$ is of type $2^+$, then
    $|\tD_g|=|\varphi_H(\tD_g)|=2$ and $|\phi_{K^+}(\tD_g)|=1$.
\end{itemize}
\end{claim}

\begin{proof}
If $g+L$ is not of type $0$ and $g\notin L$, then by \refe{cosetAD} the set
$D_g$ is disjoint with the set $A_g+A_0\seq g+L$, which contains a translate
of $A_0$. However, $A_0$ intersects non-trivially each of the four cosets of
$H$ contained in $L$. Thus, the complement of $D_g$ in $g+L$ contains an
element in each coset of $H$ contained in $g+L$. This shows the existence of
$\tD_g\seq g+L$ with $D_g\seq\tD_g$ and $|\tD_g|=|\phi_H(\tD_g)|$, and thus
proves (i).

Now suppose that $g+L$ is of type $2^-$. We assume without loss of generality
that $g\in A$ and, consequently, that either $A_g=\{g,g+a_3\}$ or
$A_g=\{g,g+a_1+a_2\}$ holds. By \refe{cosetAD}, the set $D_g$ is contained in
the complement of $A_g+A_0$ in $g+L$, which in the former case is
$\{g+a_1+a_2,g+a_1+a_2+a_3\}$, and in the latter case
$\{g+a_1+a_3,g+a_2+a_3\}$. To prove (ii), it remains to observe that each of
these sets is contained in a $K^-$-coset, but not contained in an $H$-coset.

The proof of (iii) goes along similar lines.
\end{proof}

\section{Using Two Isolated Edges: Completion of the Proof.}
  \label{s:2edgescompletion}

In this section, we complete the proof of Theorem \reft{round}. We keep the
notation and assumptions of the previous section, and since
$\frac{11}{36}\,2^r+3>\frac13\,2^r+2$ for $r\in[1,5]$, we may and do assume,
in view of Theorem \reft{roundlight}, that $r\ge 6$. To argue by
contradiction, we also assume that $|A|>\frac{11}{36}\,2^r+3$. Our goal is to
show that these assumptions are inconsistent.

\begin{claim}\label{m:caseII}
We have $\min\{n_4^-,n_4^+\}<n_0+3$.
\end{claim}

\begin{proof}
Suppose by contradiction that $n_4^-\ge n_0+3$ and $n_4^+\ge n_0+3$, and let
$A^-$ denote the union of all sets $A_g$ such that $g+L$ is of type $4^-$.
Since $|\phi_L(A)|=2^{r-3}-n_0$ and $|\phi_L(A^-)|=n_4^-\ge n_0+3$, by Lemma
\refl{php-result} every element of $\Fr/L$ is representable in at least three
ways as a sum of an element from $\phi_L(A)$ and an element from
$\phi_L(A^-)$. Hence, observing that $A^-$ is a union of $K^-$-cosets and
that each $L$-coset is a union of two $K^-$-cosets, we conclude that every
$L$-coset contains a $K^-$-coset disjoint from $D(A)$. Similarly, every
$L$-coset contains a $K^+$-coset disjoint with $D(A)$. As the union of a
$K^-$-coset and a $K^+$-coset contained in the same $L$-coset covers all this
$L$-coset with the exception of an $H$-coset, applying Claim~\refm{DinXY} we
conclude that $|D_g|\le 1$ if $g+L$ is of type
$1,2^-$, or $2^+$, and that $|D_g|\le 2$ if $g+L$ is of type $0$. Combining
this observation with Claim~\refm{Dg} and using \refe{sum-n_i-is-m/8} and
\refe{sum-n_i-is-A}, we derive
\begin{align*}
  |D(A)| &\le 2n_0+n_1+n_2^-+n_2^++2 \\
         &\le 2(n_0+n_1+n_2+n_3+n_4) - \frac12\,(n_1+2n_2+3n_3+4n_4) + 2 \\
         &=   2^{r-2} - \frac12\,|A| + 2.
\end{align*}
Compared with Corollary \refc{DishalfA}, this yields $|A|\le 2^{r-2}+2$, a
contradiction.
\end{proof}

Being the only place where the factor $\frac{11}{36}$ emerges, the following
claim can be considered the bottleneck of our method.
\begin{claim}\label{m:Case3}
We have $\max\{n_4^-,n_4^+\}<n_0+3$.
\end{claim}

\begin{proof}
Switching the notation, if necessary, and in view of Claim \refm{caseII}, we
assume by contradiction that
\begin{equation}\label{e:n4n0n4}
  n_4^-<n_0+3\le n_4^+.
\end{equation}

Let $A^+$ be the union of all sets $A_g$ such that $g+L$ is of type $4^+$. As
in the proof of Claim \refm{caseII}, every element of $\Fr/L$ is
representable in at least three ways as a sum of an element from $\phi_L(A)$
and an element from $\phi_L(A^+)$, and $A^+$ is a union of $K^+$-cosets;
hence every $L$-coset contains a $K^+$-coset disjoint from $D(A)$.
Consequently, in view of Claim \refm{DinXY}~(ii), we have $|D_g|\le 1$
whenever $g+L$ is of type $2^-$. Thus, by Claim \refm{Dg},
\begin{equation}\label{e:D-dan}
  |D(A)|\le 4n_0 + 2n_1 + n_2^- + 2n_2^+ + 2.
\end{equation}

Let $B$ denote the set of all those elements of $A$ adjacent in $\Gam(A)$ to
an element from $A^+$. As $A^+$ is a union of $K^+$-cosets, for any $b\in B$
we have $|(b+K^+)\cap A|=1$ (else $b$ could not be adjacent to an element
from $A^+$); it follows that $B$ is disjoint with $A^+$ and, since there are
precisely $n_1+2n^-_2+n_3^-$ elements $b\in A$ such that $|(b+K^+)\cap A|=1$,
that
  $$ |B| \le n_1+2n_2^-+n_3^-. $$

Consider the subgraph $\Gam'$ of $\Gam(A)$ induced by the elements of
$A^+\cup B$. Since $B$ is a vertex cover of $\Gam'$, the matching number $t'$
of $\Gam'$ does not exceed $|B|$; hence,
\begin{equation}\label{e:t'upper}
  t' \le n_1 + 2n_2^- + n_3^-.
\end{equation}

Let $t$ be the matching number of $\Gam(A)$ and let $T$ be a matching in
$\Gam(A)$ with $|T|=t$ edges.
As $\Gamma(A)$ has no isolated vertices, the number of edges between $A^+$
and $B$ in $\Gam(A)$ is at least $|A^+|=4n_4^+$, and the definition of $t'$
ensures that at most $t'$ of these edges belong to $T$; thus,
\begin{equation}\label{e:meow}
  t \le |D(A)| - 4n_4^+ + t'.
\end{equation}

To obtain another relation between $t$ and $t'$, we notice that if $b\in B$
is adjacent in $\Gam(A)$ to $a\in A^+$, then in fact every element of
$D_{a+b}$ corresponds to an edge incident with $b$: for all elements of
$D_{a+b}$ are contained in a $K^+$-coset (as shown earlier), and since
$a+b\in D_{a+b}$, this coset is $a+b+K^+=b+A_a$. Now, fix a matching $T'$ of
$\Gam'$. As any edge of $T$ corresponds to an element from $D(A)$, we see
from Claim \refm{Dg} that corresponding to the edges of $T$ are at most four
elements from every $L$-coset of type $0$, two elements from $L$, and at most
two elements from every $L$-coset of type $1$, $2^-$ or $2^+$. Taking into
account that, for each edge $(a,b)$ of $T'$, there is actually at most one
element in the coset $a+b+L$ corresponding to an edge of $T$ (as all edges in
$D_{a+b}$ are adjacent to the same vertex), we conclude that
\begin{equation}\label{e:ruuff}
   t \le 4n_0 + 2n_1 + 2n_2^- + 2n_2^+ + 2 - t'.
\end{equation}

We complete the proof of Claim \refm{Case3} showing that an appropriate
combination of the estimates \refe{n4n0n4}--\refe{ruuff} yields a
contradiction to the assumption on the size of $A$, made at the beginning of
this section. Specifically, substituting \refe{sum-n_i-is-A} into the
estimate of Corollary~\refc{matching}, we get
  $$ n_1 + 2n_2 + 3n_3 + 4n_4 \le |D(A)| + t - 4. $$
Taking the sum of this inequality with weight $4$, identity
\refe{sum-n_i-is-m/8} with  weight $44$, the first inequality in
\refe{n4n0n4} in the form $-n_0+n_4^-\le 2$ with weight $12$, and
inequalities \refe{D-dan}, \refe{t'upper}, \refe{meow}, and \refe{ruuff} with
weights $7$, $2$, $3$ and $1$, respectively, we obtain
  $$ 30n_1 + 52n_2^0 + 39n_2^- + 36n_2^+ + 54n_3^- + 56n_3^+
                             + 72n_4^- + 72n_4^+ \le 44\cdot 2^{r-3} - 20. $$
(The weights were found by solving the corresponding linear program to yield
the best possible bound.) In view of \refe{sum-n_i-is-A}, the left-hand side
is at least as large as
  $$ 18(n_1+2n_2+3n_3+4n_4) = 18(|A|-4); $$
thus
  $$ |A| \le \frac{22}9\cdot 2^{r-3} - \frac{10}9 + 4
                                            < \frac{11}{36}\cdot 2^r + 3, $$
a contradiction.
\end{proof}

\begin{claim}\label{m:caseI}
We have $n_4\ge n_0+n_1+n_2^0+4$.
\end{claim}

\begin{proof}
Applying Claim \refm{Dg} and using \refe{sum-n_i-is-m/8} and
\refe{sum-n_i-is-A}, we get
\begin{align}
  |D(A)|
    &\le 4n_0+2n_1+2n_2^-+2n_2^++2 \nonumber \\
    &=   6(n_0+n_1+n_2+n_3+n_4)-2(n_1+2n_2+3n_3+4n_4) \nonumber \\
    &\phantom{6(n_0+n_1+n_2+n_3+n_4)-\ \ } + 2(n_4-n_0-n_1-n_2^0+1)
                                                               \nonumber \\
    &=   3\cdot 2^{r-2}-2|A| + 2(n_4-n_0-n_1-n_2^0+2). \label{e:r=6A=23}
\end{align}
Comparing with Corollary~\refc{DishalfA}, we obtain
  $$ \frac12\,|A| \le 3\cdot 2^{r-2}-2|A| + 2(n_4-n_0-n_1-n_2^0+2). $$
Hence
  $$ n_4-n_0-n_1-n_2^0 \ge \frac54|A| - 3\cdot2^{r-3} - 2, $$
and if $r\ge 7$, then the result follows from
  $$ |A|\ge \lcl \frac{11}{36}\,2^r + 3 \rcl > \frac{3}{5}\,2^{r-1} + 4. $$
If $r=6$ and $|A|\ge 24$, then the estimate $|A|>\frac{3}{5}\,2^{r-1}+4$
remains valid, proving the result in this case, too.

In view of the assumption $r\ge 6$ made at the beginning of this section, and
since for $r=6$ we have $\lcl \frac{11}{36}\,2^r+3\rcl=23$, we are left with
the case where $r=6$ and $|A|=23$, which we proceed to consider. By
Corollary~\refc{DishalfA}, we have $|D(A)|\ge 12$, whence
Proposition~\refp{line} gives
  $$ \deg(a_1)+\deg(a_2) \ge |A|+|D(A)|-2^{r-1} \ge 3 $$
for every edge $(a_1,a_2)$ of $\Gam(A)$. Now Lemma \refl{linegraph}, applied
with $\del=3$, yields $|D(A)|\ge\frac23\,|A|$, and hence in fact
 $|D(A)|\ge 16$. Substituting into \refe{r=6A=23}, we get
  $$ 16 \le 3\cdot16 - 2\cdot 23 + 2 (n_4-n_0-n_1-n_2^0+2), $$
leading to
  $$ n_4-n_0-n_1-n_2^0+2 \ge 7 $$
and implying the result.
\end{proof}

\begin{claim}\label{m:matching}
The matching number of $\Gam(A)$ does not exceed $n_1+2n_2^-+2n_2^++n_3+2$.
\end{claim}

\begin{proof}
Write $\sig:=a_1+a_2+a_3$, so that $H=\{0,\sig\}$. Observe that
 $\sig\notin D(A)$, in view of Proposition \refp{tfree} and since
$a_1,a_2+a_3\in D(A)$. If $(a,b)$ is an edge in $\Gam(A)$, then either
$a+\sig\notin A$ or $b+\sig\notin A$: otherwise
$a+b\in D(A)$ would be represented as $(a+\sig)+(b+\sig)$ with both summands
in $A$ and distinct from $a$ and $b$. This shows that every edge of $\Gam(A)$
is incident with an element of the set $B:=\{a\in A\colon a+\sig\notin A\}$.
However, by the
definition of the quantities $n_i,n_i^0,n_i^+$, and $n_i^-$, the total number
of elements of $B$ is
  $$ n_1 + 2n_2^- + 2n_2^+ + n_3 + 4. $$
It remains to notice that, in a matching of $\Gam(A)$, no two
distinct edges can be incident to the same element of $B$, and that any
maximal matching contains the isolated edges $(0,a_1)$ and $(a_2,a_3)$,  both
incident to two elements of $B$.
\end{proof}

By Claims \refm{Case3} and \refm{caseI}, it remains to consider the case
where
\begin{equation}\label{e:n4-+n0}
  n_4^-\le n_0+2,\ n_4^+\le n_0+2
\end{equation}
and
\begin{equation}\label{e:claim6}
  n_4 \ge n_0 + n_1 + 4,
\end{equation}
which from now on we assume to hold. Notice that these assumptions imply
\begin{equation}\label{e:atleasttwo}
  \min\{n_4^-, n_4^+\}\ge 2.
\end{equation}
As above, we define $A^-$ to be the union of all sets $A_g$ such that $g+L$
is of type $4^-$, and $A^+$ to be the union of those $A_g$ with $g+L$ of type
$4^+$. Next, let $B$ be the union of all $A_g$ with $|A_g|\ge 2$; thus,
\begin{equation}\label{e:phiL}
  |\phi_L(A^-)|=n_4^-,\ |\phi_L(A^+)|=n_4^+,
                                \ \text{and}\ |\phi_L(B)|=n_2+n_3+n_4+1.
\end{equation}
Furthermore, let $C^-$ denote the set of all those $g\in A^-+B$ with the
property that $\phi_L(g)$ has at least two representations as an element from
$\phi_L(A^-)$ and an element from $\phi_L(B)$; similarly, denote by $C^+$ the
set of those $g\in A^++B$ such that $\phi_L(g)$ has at least two
representations as an element from $\phi_L(A^+)$ and an element from
$\phi_L(B)$.

If $A_g\seq A^-$, for some $g\in\Fr$, and if $b_1,b_2\in B$ are distinct and
belong to the same $L$-coset, then the $K^-$-cosets $b_1+A_g$ and $b_2+A_g$
either coincide or cover the whole $L$-coset $g+\{b_1,b_2\}+L$. It follows
that if $g\in C^-+L$, then $g+L$ contains a $K^-$-coset disjoint from $D(A)$.
Likewise, if $g\in C^++L$, then $g+L$ contains a $K^+$-coset disjoint from
$D(A)$. Hence, for $g\in(C^-+L)\cap(C^++L)$, the set $D_g$ is contained in an
$H$-coset, and thus, by Claim~\refm{DinXY},
\begin{equation}\label{e:dust}
  |D_g| \le \begin{cases}
              1 &\text{if $g+L$ is of type $1$ or $2$,} \\
              2 &\text{if $g+L$ is of type $0$.}
            \end{cases}
\end{equation}

By the pigeonhole principle, we have
\begin{equation}\label{e:F-F+box}
  |\phi_L(C^-)\cap\phi_L(C^+)| \ge |\phi_L(C^-)|+|\phi_L(C^+)|-2^{r-3},
\end{equation}
while, by Corollary \refc{S2}
and \refe{phiL},
\begin{align}
  |\phi_L(C^-)| &\ge \min \{ 2n_2+2n_3+2n_4+2n_4^--2^{r-3}-2, n_2+n_3+n_4 \}
                                                                    \nonumber
\intertext{and}
  |\phi_L(C^+)| &\ge \min \{ 2n_2+2n_3+2n_4+2n_4^+-2^{r-3}-2, n_2+n_3+n_4 \}.
                                                          \label{e:secondmin}
\end{align}

We notice that at least one of these minima is attained on its second term,
for if
\begin{align*}
  2n_2+2n_3+2n_4+2n_4^--2^{r-3}-2 \le n_2+n_3+n_4
\intertext{and}
  2n_2+2n_3+2n_4+2n_4^+-2^{r-3}-2 \le n_2+n_3+n_4
\end{align*}
both hold true, then taking their sum we obtain
  $$ 2n_2 + 2n_3 + 4n_4 \le 2^{r-2}+4,$$
which, in view of \refe{sum-n_i-is-m/8}, can be re-written as
  $$ n_4 \le n_0 + n_1 + 3; $$
this, however, is inconsistent with \refe{claim6}.

By symmetry, we can assume that
\begin{equation}\label{e:phiFminus}
  |\phi_L(C^-)| \ge n_2 + n_3 + n_4,
\end{equation}
and we consider two cases, according to the value in the right-hand side of
\refe{secondmin}.

If $|\phi_L(C^+)|\ge n_2+n_3+n_4$, then
from \refe{F-F+box} and \refe{phiFminus} we derive
  $$ |\phi_L(C^-)\cap\phi_L(C^+)| \ge 2n_2+2n_3+2n_4-2^{r-3}. $$
Consequently, there are at least
  $$ (2n_2+2n_3+2n_4-2^{r-3}) - (n_3+n_4+1) = 2n_2+n_3+n_4-2^{r-3}-1 $$
$L$-cosets of type $0,1$, or $2$ contained in $(C^-+L)\cap(C^++L)$. Hence, by
Claim~\refm{Dg} and
estimate \refe{dust}, we see that
\begin{align*}
  |D(A)| &\le 4n_0+2n_1+2n_2+2 - (2n_2+n_3+n_4-2^{r-3}-1) \\
         &=   4n_0 + 2n_1 - n_3 - n_4 + 2^{r-3} + 3.
\end{align*}
Combining this estimate with Corollary~\refc{matching} and
Claim~\refm{matching}, we get
  $$ |A| \le 4n_0 + 3n_1 + 2n_2^- + 2n_2^+ - n_4 + 2^{r-3} + 5 $$
and furthermore, substituting the value of $|A|$ from \refe{sum-n_i-is-A},
  $$ -4n_0 - 2n_1 + 2n_2^0 + 3n_3 + 5n_4 \le 2^{r-3} + 1. $$
Taking the sum of this
estimate, inequalities \refe{n4-+n0}, and identity \refe{sum-n_i-is-m/8}
multiplied by $6$, we obtain
  $$ 4n_1 + 2n_2^0 + 6n_2 + 9n_3 + 12n_4 \le 7\cdot 2^{r-3} - 1. $$
By \refe{sum-n_i-is-A}, the expression in the left-hand side is at least
$3(|A|-4)$; consequently,
  $$ |A| \le \frac7{24}\,2^r - \frac13 + 4 < \frac{11}{36}\,2^r + 3 $$
(in view of $r\ge 6$), a contradiction.

Finally, suppose that $|\phi_L(C^+)|\ge 2n_2+2n_3+2n_4+2n_4^+-2^{r-3}-2$.
Arguing as in the previous case, we get
\begin{align*}
  |\phi_L(C^-)\cap\phi_L(C^+)|
    &\ge 3n_2 + 3n_3 + 3n_4 + 2n_4^+ - 2^{r-2} - 2, \\
  |D(A)|
    &\le 4n_0 + 2n_1 + 2n_2 + 2 \\
    &\phantom{4n_0+}
              - \big((3n_2+3n_3+3n_4+2n_4^+-2^{r-2}-2)-(n_3+n_4+1)\big) \\
    &=   4n_0 + 2n_1 - n_2 - 2n_3 - 2n_4 - 2n_4^+ + 2^{r-2} + 5, \\
  |A|
    &\le 4n_0 + 3n_1 - n_2^0 + n_2^- + n_2^+ - n_3 - 2n_4
                                              - 2n_4^+ + 2^{r-2} + 7,
\end{align*}
and hence
  $$ -4n_0 - 2n_1 + 3n_2^0 + n_2^- + n_2^+ + 4n_3 + 6n_4 + 2n_4^+
                                                        \le 2^{r-2} + 3. $$
Taking the sum of the last inequality, the first of the inequalities
\refe{n4-+n0}, and identity \refe{sum-n_i-is-m/8} multiplied by $5$, we
obtain
  $$ 3n_1 + 2n_2^0 + 6n_2 + 9n_3 + 12n_4 + n_4^+ \le 7\cdot 2^{r-3}. $$
In view of \refe{atleasttwo} and \refe{sum-n_i-is-A}, this yields
  $$ 3(|A|-4) \le 7\cdot 2^{r-3} - 2, $$
leading to a contradiction as above. This completes the proof of
Theorem~\reft{round}.

\section*{Acknowledgement}

We are grateful to Alexander Davydov for attracting our attention to the
problem of studying $1$-saturating sets and for several useful remarks,
including mentioning to us the exceptional minimal $1$-saturating set of size
$11$ in $\F_2^5$.

\bigskip

\end{document}